\documentclass[draft,12pt]{amsart}
\usepackage{amssymb}

\theoremstyle{plain}
\newtheorem{lemma}{Lemma}[section]

\newtheorem{theorem}[lemma]{Theorem}
\newtheorem{corollary}[lemma]{Corollary}
\newtheorem{conjecture}[lemma]{Conjecture}

\theoremstyle{definition}

\newtheorem{remark}[lemma]{Remark}

\newcommand{\ggcaffil}[1]{\dedicatory{\textup{\larger{#1}}}}

\newcommand{\chl}{\ensuremath{\chi_{_\ell}}}
\newcommand{\lm}{\ensuremath{\lambda}}

\newcommand{\ggcqedsymbol}{$\square$}
\newcommand{\ggcqed}{\hbox{}\nobreak\hbox{\quad\ggcqedsymbol}}
\newcommand{\ggcendpf}{\ggcqed}
\newcommand{\ggcnopf}{\ggcqed}
\newcommand{\ggcendconj}{\ggcqed}
\newcommand{\ggcenddef}{\ggcqed}

\newcommand{\ggcendrem}{\ggcenddef}

\begin{document}
\title[A Lower Bound for Partial List Colorings]
  {A Lower Bound for Partial List Colorings}
\author{Glenn G. Chappell}
\ggcaffil{Department of Mathematics, Southeast Missouri
  State University}
\address{Department of Mathematics\\
  Southeast Missouri State University\\
  Cape Girardeau, MO 63701\\
  USA}
\email{gchappell@semovm.semo.edu}
\subjclass{05C15}
\date{May 13, 1998}
\begin{abstract}
  Let $G$ be an $n$-vertex graph with list-chromatic number $\chl$.
  Suppose each vertex of $G$ is assigned a list of $t$ colors.
  Albertson, Grossman, and Haas~\cite{AGH98} conjecture
  that at least $\frac{t n}{\chl}$ vertices can be colored from
  these lists.
  We prove a lower bound for the number of colorable vertices.
  As a corollary, we show that at least $\frac{6}{7}$
  of the conjectured number can be colored.
\end{abstract}
\maketitle

\section{Introduction} \label{S:intro}

Let $G$ be an $n$-vertex graph.
Let each vertex of $G$ be assigned a list of $s$ colors.
An \emph{$s$-list coloring} of $G$ is a proper vertex coloring
in which each vertex is given a color from its list.
We say $G$ is \emph{$s$-choosable} if $G$ has an $s$-list coloring
for each assignment of lists of size $s$ to the vertices.
The \emph{list-chromatic number} of $G$,
denoted by $\chl$,
is the least $s$ so that $G$ is $s$-choosable.
List coloring was introduced independently
by Vizing~\cite{VizV76}
and by Erd\H os, Rubin, and Taylor~\cite{ERT79}.

Suppose each vertex of $G$ is assigned
a list of $t$ colors ($t\le\chl$).
It may not be possible to color
every vertex from these lists;
however, we can color a subset of the vertices.
Albertson, Grossman, and Haas~\cite{AGH98} ask
how many vertices can be colored.
Letting $\ell_t$ range over all assignments
of lists of size $t$ to the vertices of $G$,
we define
\[
\lm_t:=\min_{\ell_t}
  \left\{
  \text{max.~no.~of vertices of $G$ colorable from the
    lists $\ell_t$}
  \right\}.
\]

\begin{conjecture}[Albertson, Grossman, \& Haas {\cite{AGH98}}]
  \label{J:aghmain}
If $t\le\chl$, then
$\lm_t\ge\frac{t n}{\chl}$.\ggcendconj\end{conjecture}

Conjecture~\ref{J:aghmain} clearly holds when $t=0$ or $t=\chl$.
When $t=1$, we note that $G$ has an independent set of size
at least $\frac{n}{\chi}\ge\frac{n}{\chl}$, where $\chi$ is
the ordinary chromatic number.
Since an independent set can always be colored from lists
of size 1, we have $\lm_1\ge\frac{n}{\chl}$, and so
the conjecture holds for $t=1$.
However, the conjecture remains open for $1<t<\chl$.

We will use a method similar to that of the proof
of~\cite[Thm.~2]{AGH98} to prove a lower bound for $\lm_t$.
As a corollary, we show that
$\lm_t>\frac{6}{7}\cdot\frac{t n}{\chl}$.

\section{The Results} \label{S:results}

For integers $s>t>0$, we define
\[
f_{s,t}(x):=1-x-\left[1-\frac{1-x}{s-t}\right]^t.
\]
Then $f_{s,t}(0)>0$, $f_{s,t}(1)<0$,
and $\frac{d}{dx}f_{s,t}(x)<0$ for $x>0$.
Thus, the equation $f_{s,t}(x)=0$ has exactly one positive
solution, which we denote by $q_{s,t}$.
Note that $0<q_{s,t}<1$.

\begin{theorem} \label{T:lowerbd}
If $G$ is $s$-choosable, and $s>t>0$,
then $\lm_t\ge q_{s,t}\cdot n$.
Furthermore, this inequality is strict if $t>1$.
\end{theorem}

\begin{proof}
Let $q=q_{s,t}$.
Let each vertex of $G$ be assigned a list of $t$ colors,
and let $R$ be the union of these lists.
Let $u=s-t$.
We augment each list by adding new colors
$\pi_1,\pi_2,\dots,\pi_u\not\in R$.
There is an $s$-list coloring $\varphi$ from the augmented lists.
For $1\le i\le u$, let $I_i$ be the set of
vertices given color $\pi_i$;
each $I_i$ is an independent set.
Let $\mathcal I=I_1\cup\cdots\cup I_u$,
and let $H=G-\mathcal I$.
Thus, $\varphi$ colors the vertices of $H$ with colors in $R$.

We partition the set $R$ randomly into classes
$R_0,R_1,\dotsc,R_u$;
the class in which each color falls is chosen independently.
We place a given color in class $R_0$ with probability $q$ and in each
of the other classes with probability $\frac{1-q}{u}$.

Given a partition of $R$, we construct a proper
$t$-list coloring of a subset
of the vertex set of $G$ as follows.
The colors in $R_0$ are used on the vertices of $H$.
We color a vertex in $H$ with its color in the coloring $\varphi$ if
that color appears in $R_0$;
otherwise, we leave the vertex uncolored.
Thus, each vertex in $H$ is colored with probability $q$.

For $1\le i\le u$, the colors in $R_i$
are used on $I_i$.
We color a vertex in $I_i$ if a color in its list appears in $R_i$;
otherwise, we leave the vertex uncolored.
A particular color fails to appear in $R_i$ with
probability
\[
1-\frac{1-q}{u}.
\]
Thus, each vertex of $\mathcal I$ is colored with
probability
\[
1-\left[1-\frac{1-q}{u}\right]^t
=f_{s,t}(q)+q=q,
\]
since $f_{s,t}(q)=0$.

For some partition of $R$
the number of colored vertices
is at least the expected number.
Since the above scheme results in each vertex being colored
with probability $q$,
we have $\lm_t\ge qn$.

Now let $t>1$.
To see that the above inequality is strict,
note that $q$ is a root of the polynomial
$p(x):=u^t f_{s,t}(x)$.
All the coefficients of $p(x)$ are integers;
the leading coefficient of is $-1$.
By the Rational Roots Theorem, every rational root of $p(x)$
is an integer.
Since $0<q<1$,
we conclude that $q$ is irrational, and so
$\lm_t\ne qn$.\ggcendpf\end{proof}

Albertson et al.\ showed~\cite[Thm.~2]{AGH98}
that $\lm_2>\frac{3n}{5}$ if $\chl=3$.
Theorem~\ref{T:lowerbd} allows us to improve on this.
Since $q_{3,2}=\frac{-1+\sqrt{5}}{2}$,
we have the following stronger result.

\begin{corollary} \label{C:3,2}
If $\chl=3$, then
$\lm_2>\frac{-1+\sqrt{5}}{2}\cdot n$.\ggcnopf
\end{corollary}

The conjectured lower bound is greater still:
$\frac{2n}{3}$.

\begin{remark} \label{R:planar5,4}
Thomassen~\cite{ThoC94} showed that
every planar graph is $5$-choosable.
Thus, we can apply Theorem~\ref{T:lowerbd} to planar
graphs.
Albertson et al.~\cite{AGH98} estimate
that their method could be used to
show that $\lm_4\ge\frac{7}{10}n$ for every planar graph.
By definition, $q_{5,4}$ is the positive root
of $1-x-x^4$: a bit more than $0.724$.
Thus, $\lm_4>0.724n$ if $\chl\le5$ and
in particular for every planar graph.\ggcendrem\end{remark}

Next we show that $\lm_t$
is at least $\frac{6}{7}$
of the conjectured lower bound.

\begin{corollary} \label{C:nicelowerbd}
If $t\le\chl$, then
$\lm_t>\frac{6}{7}\cdot\frac{t n}{\chl}$.\ggcnopf\end{corollary}

Corollary~\ref{C:nicelowerbd} follows
from Theorem~\ref{T:lowerbd},
the comments after Conjecture~\ref{J:aghmain},
and the following lemma.

\begin{lemma} \label{L:whatisqst}
If $s>t>0$, then
\(
\frac{6}{7}\cdot\frac{t}{s}
  <q_{s,t}\le\frac{t}{s}.
\)
\end{lemma}

\begin{proof}
We first show that
\[
\frac{6}{7}\cdot\frac{t}{s}<q_{s,t}.
\]
Recall that $f_{s,t}(q_{s,t})=0$,
and that $f_{s,t}(x)$ is decreasing for $x>0$.
Thus, it suffices to show that
$f_{s,t}(cv)>0$,
where $c=\frac{6}{7}$ and $v=\frac{t}{s}$.
\begin{align*}
f_{s,t}\left(c v\right)
  &=1-c v-\left[1-\frac{1-c v}{s-s v}\right]^{s v}\\
&>1-c v-e^{-v\cdot\frac{1-c v}{1-v}},\quad
  \text{since $\left(1-\frac{a}{b}\right)^b<e^{-a}$, for $0<a<b$.}
\end{align*}
We show that the last expression is positive.
Note that both
$1-c v$ and $e^{-v\cdot\frac{1-c v}{1-v}}$ are positive;
let $g(v)$ be the difference of the logarithms of these quantities.
\[
g(v) := \ln\left(1-c v\right)-\left(-v\cdot\frac{1-c v}{1-v}\right).
\]
It suffices to show that $g(v)>0$, for $0<v<1$.
Verification of this requires only the standard techniques
of calculus and is omitted.

Now we show that
\[
q_{s,t}\le\frac{t}{s}.
\]
Since $f_{s,t}(x)$ is decreasing for $x>0$,
it suffices to show that $f_{s,t}\left(\frac{t}{s}\right)\le0$.
\begin{align*}
f_{s,t}\left(\frac{t}{s}\right)
  &=1-\frac{t}{s}-\left[1-\frac{1-\frac{t}{s}}{s-t}\right]^t\\
&=1-\frac{1}{s}\cdot t-\left[1-\frac{1}{s}\right]^t\\
&\le0.\ggcendpf
\end{align*}
\end{proof}

The techniques of the above proof can be used
to prove a slightly better lower bound for $q_{s,t}$.
Numerical computations indicate that
\[
\inf_{s,t}\frac{q_{s,t}}{t/s}\approx0.8598841287>\frac{6}{7}.
\]

\end{document}